\documentclass[11pt]{article}%
\usepackage{amsmath}
\usepackage{amsfonts}
\usepackage{amssymb}
\usepackage{graphicx}%
\makeatletter
\newfont{\footsc}{cmcsc10 at 8truept}
\newfont{\footbf}{cmbx10 at 8truept}
\newfont{\footrm}{cmr10 at 10truept}
\newtheorem{theorem}{Theorem}

\newtheorem{corollary}[theorem]{Corollary}

\newtheorem{proposition}[theorem]{Proposition}

\newenvironment{proof}[1][Proof]{\noindent{\textbf {#1}  }}  {\hfill$\square$\bigskip}

\def\blfootnote{\xdef\@thefnmark{}\@footnotetext}
\begin{document}

\title{Maxima of the $Q$-index: degenerate graphs}
\author{Vladimir Nikiforov\thanks{Department of Mathematical Sciences, University of
Memphis, Memphis TN 38152, USA; \textit{email: vnikifrv@memphis.edu}}}
\maketitle

\begin{abstract}
Let $G$ be a $k$-degenerate graph of order $n.$ It is well-known that $G\ $has
no more edges than $S_{n,k},$ the join of a complete graph of order $k$ and an
independent set of order $n-k.$ In this note it is shown that $S_{n,k}$ is
extremal for some spectral parameters of $G$ as well. More precisely, letting
$\mu\left(  H\right)  $ and $q\left(  H\right)  $ denote the largest
eigenvalues of the adjacency matrix and the signless Laplacian of a graph $H,$
the inequalities%
\[
\mu\left(  G\right)  <\mu\left(  S_{n,k}\right)  \text{ \ \ \ and
\ \ \ }q\left(  G\right)  <q\left(  S_{n,k}\right)
\]
hold, unless $G=S_{n,k}$.

The latter inequality is deduced from the following general bound, which
improves some previous bounds on $q\left(  G\right)  $:

If $G$ is a graph of order $n$, with $m$ edges, with maximum degree $\Delta$
and minimum degree $\delta,$ then%
\[
q\left(  G\right)  \leq\min\left\{  2\Delta,\frac{1}{2}\left(  \Delta
+2\delta-1+\sqrt{\left(  \Delta+2\delta-1\right)  ^{2}+16m-8\left(
n-1+\Delta\right)  \delta}\right)  \right\}  .
\]
Equality holds if and only if $G$ is regular or $G$ has a component of order
$\Delta+1$ in which every vertex is of degree $\delta$ or $\Delta,$ and all
other components are $\delta$-regular.\medskip



\end{abstract}

\section{Introduction and main results}

A graph $G$ is said to be\emph{ }$k$\emph{-degenerate} if every subgraph of
$G$ contains vertex of degree at most $k.$ Clearly, $k$-degenerate graphs
constitute a monotone graph family, as every subgraph of a $k$-degenerate
graph is $k$-degenerate itself. In fact, many much studied classes of graphs
consist of degenerate graphs: e.g., forests are $1$-degenerate, planar graphs
are $5$-degenerate, graphs of tree-width $k$ are $k$-degenerate, etc.

By induction on $n$ one can show that if $n\geq k$ and $G$ is a $k$-degenerate
graph of order $n,$ then the number of edges $e\left(  G\right)  $ of $G$
satisfies%
\[
e\left(  G\right)  \leq kn-\left(  k^{2}+k\right)  /2.
\]
Equality holds for numerous graphs, in particular if $G=S_{n,k},$ where
$S_{n,k}$ is the join of a complete graph of order $k$ and an independent set
of order $n-k.$

It turns out that the graph $S_{n,k}$ is extremal for some spectral parameters
of $k$-degenerate graphs as well. A notable early contribution is due to Hong,
who showed in \cite{Hon05} that if $G$ has tree-width $k,$ then the largest
eigenvalue $\mu\left(  G\right)  $ of the adjacency matrix of $G$ satisfies
$\mu\left(  G\right)  <\mu\left(  S_{n,k}\right)  ,$ unless $G=S_{n,k}.$ Here
we prove the following stronger theorem.

\begin{theorem}
\label{tA}If $G$ is a $k$-degenerate graph of order $n\geq k$, then
$\mu\left(  G\right)  <\mu\left(  S_{n,k}\right)  ,$ unless $G=S_{n,k}.$
\end{theorem}

We shall show that Theorem \ref{tA} readily follows from a well-known bound on
$\mu\left(  G\right)  ,$ stated below as Theorem A. The main effort of this
note, however, is to prove a similar bound on the largest eigenvalue $q\left(
G\right)  $ of the signless Laplacian of a $k$-degenerate graph $G.$

\begin{theorem}
\label{tQ}If $G$ is a $k$-degenerate graph of order $n\geq k$, then $q\left(
G\right)  <q\left(  S_{n,k}\right)  ,$ unless $G=S_{n,k}.$
\end{theorem}

The simple proof strategy for Theorem \ref{tA} did not apply to Theorem
\ref{tQ} because until recently there was no appropriate bound on $q\left(
G\right)  $ similar to Theorem A. To fill in this void we came up with the
following theorem, which is of separate interest and may have other
applications as well.

\begin{theorem}
\label{tI}If $G$ is a graph of order $n$, with $m$ edges, with maximum degree
$\Delta$ and minimum degree $\delta,$ then
\begin{equation}
q\left(  G\right)  \leq\min\left\{  2\Delta,\frac{\Delta+2\delta
-1+\sqrt{\left(  \Delta+2\delta-1\right)  ^{2}+16m-8\left(  n-1+\Delta\right)
\delta}}{2}\right\}  .\label{bo}%
\end{equation}
Equality holds if and only if $G$ is regular or $G$ has a component of order
$\Delta+1$ in which every vertex is of degree $\delta$ or $\Delta,$ and all
other components are $\delta$-regular.
\end{theorem}

Inequality (\ref{bo}) seems quite unwieldy at first glance, so some remarks
are due here. First, obviously the right side is nondecreasing in $m;$ a
simple calculation shows that it is equal to $2\Delta$ whenever%
\[
2m\geq\Delta^{2}+\Delta+\left(  n-1-\Delta\right)  \delta\text{.}%
\]

Furthermore, the dependence of the right side of (\ref{bo}) on $\delta$ and
$\Delta$ is clarified by the following proposition whose simple proof is omitted.

\begin{proposition}
\label{pro1}Let $n\geq1$ and $0\leq m\leq n\left(  n-1\right)  /2.$ If
$x\leq2m/n$ and $y\geq2m/n,$ then the functions
\[
f\left(  x,y\right)  =y+2x-1+\sqrt{\left(  y+2x-1\right)  ^{2}+16m-8\left(
n-1+y\right)  x}%
\]
and
\[
g\left(  x,y\right)  =\min\left\{  2y,\frac{1}{2}\left(  y+2x-1+\sqrt{\left(
y+2x-1\right)  ^{2}+16m-8\left(  n-1+y\right)  x}\right)  \right\}
\]
are decreasing in $x$ and increasing in $y.$
\end{proposition}

Now, using Proposition \ref{pro1} and the fact that $\Delta\leq n-1,$ we
obtain the following simplified version of Theorem \ref{tI}.

\begin{corollary}
\label{cor1}If $G$ is a graph of order $n$, with $m$ edges and with minimum
degree $\delta,$ then$.$
\[
q\left(  G\right)  \leq\frac{1}{2}\left(  n+2\delta-2+\sqrt{\left(
n+2\delta-2\right)  ^{2}+16m-16\left(  n-1\right)  \delta}\right)  .
\]
Equality holds if and only if each degree of $G$ is equal to $\delta$ or
$n-1.$
\end{corollary}

Likewise, using that $\delta\geq0,$ we obtain another bound.

\begin{corollary}
If $G$ is a graph of order $n$, with $m$ edges and with maximum degree
$\Delta,$ then%
\[
q\left(  G\right)  \leq\min\left\{  2\Delta,\frac{1}{2}\left(  \Delta
-1+\sqrt{\left(  \Delta-1\right)  ^{2}+16m}\right)  \right\}  .
\]
Equality holds if and only if $G$ is $\Delta$-regular with possibly some
isolated vertices.
\end{corollary}

Similar bounds on $q\left(  G\right)  $ have been known in the literature, and
it seems that Theorem \ref{tI} compares favorably with most of them. For
example, Theorem \ref{tI} improves the bound
\[
q\left(  G\right)  \leq\frac{1}{2}\left(  \delta-1+\sqrt{\left(
\delta-1\right)  ^{2}+16m+8\Delta^{2}-8\left(  n-1\right)  \delta}\right)  ,
\]
obtained in \cite{LiPa04}, and the bound
\[
q\left(  G\right)  \leq\frac{1}{2}\left(  \Delta+\delta-1+\sqrt{\left(
\Delta+\delta-1\right)  ^{2}+16m-8\left(  n-1\right)  \delta}\right)  ,
\]
obtained in \cite{LLT04}.

Very likely Theorem \ref{tI} can be used for other extremal problems about
$q\left(  G\right)  ,$ just like Theorem A below has been used for many
extremal problems about $\mu\left(  G\right)  .$

The rest of the paper is dedicated to the proofs of Theorems \ref{tA},
\ref{tQ} and \ref{tI}.

\section{Proofs}

For graph notation and concepts undefined here, we refer the reader to
\cite{Bol98}. For introductory material on the signless Laplacian see the
survey of Cvetkovi\'{c} \cite{C10} and its references. In particular, for a
graph $G$ we write:

- $V\left(  G\right)  $ and $E\left(  G\right)  $ for the sets of vertices and
edges of $G;$

- $v\left(  G\right)  $ for $\left\vert V\left(  G\right)  \right\vert $ and
$e\left(  G\right)  $ for $\left\vert E\left(  G\right)  \right\vert $;

- $d_{u}$ for the degree of a vertex $u.\medskip$

In the proof of Theorem \ref{tA} we shall use the following inequality proved
in \cite{Nik02}, with condition for  equality proved in \cite{ZhCh05}.\medskip

\textbf{Theorem A }\emph{If }$G$\emph{ is a graph of order }$n$\emph{, with
}$m$\emph{ edges, and with minimum degree }$\delta,$\emph{ then}$.$\emph{ }%
\begin{equation}
q\left(  G\right)  \leq\frac{\delta-1}{2}+\sqrt{2m-n\delta+\frac{\left(
\delta+1\right)  }{4}^{2}}.\label{in}%
\end{equation}
\emph{Equality holds if and only if }$G$\emph{ is regular or }$G$\emph{ has a
component of order }$\Delta+1$\emph{ in which every vertex is of degree
}$\delta$\emph{ or }$\Delta,$\emph{ and all other components are }$\delta
$\emph{-regular.\medskip}

For connected graphs inequality (\ref{in}) has been proved independently by
Hong, Shu and Fang\ in \cite{HSF01}, however, here the unnecessary requirement
for connectedness blurs the situation, as equality holds in (\ref{in}) for a
much smaller set of graphs.\medskip

We shall also need the following characterization of $k$-degenerate
graphs:\medskip

\textbf{Lemma B }\emph{A graph }$G$\emph{ is }$k$\emph{-degenerate if and only
if it vertices can be ordered so that every vertex has at most }$k$\emph{
preceding neighbors.}\medskip

Finally, note the following facts about $S_{n,k}$:%
\begin{align}
\mu\left(  S_{n,k}\right)   &  =\frac{k-1}{2}+\sqrt{kn-\frac{3k^{2}+2k-1}{4}%
};\label{Smu}\\
q\left(  S_{n,k}\right)   &  =\frac{1}{2}\left(  n+2k-2+\sqrt{\left(
n+2k-2\right)  ^{2}-8\left(  k^{2}-k\right)  }\right)  .\label{Squ}%
\end{align}
\medskip\textbf{ }

\begin{proof}
[\textbf{Proof of Theorem \ref{tA}}]Our first goal is to prove that
$\mu\left(  G\right)  \leq\mu\left(  S_{n,k}\right)  .$ Without loss of
generality we may assume that $G$ is edge maximal, that is to say, if $H$ is
obtained by adding a new edge to $G,$ then $H$ is not $k$-degenerate. Order
the vertices of $G$ as follows: choose a vertex $v\in V\left(  G\right)  $
with $d\left(  v\right)  =\delta\left(  G\right)  $ and set $v_{n}=v;$
further, having chosen $v_{n},\ldots,v_{j+1},$ let $H$ be the graph induced by
the set $V\left(  G\right)  \backslash\left\{  v_{n},\ldots,v_{j+1}\right\}
,$ choose $v\in V\left(  H\right)  $ with $d_{H}\left(  v\right)
=\delta\left(  H\right)  $ and set $v_{j}=v.$

Since $G$ is $k$-degenerate, in the ordering $v_{n},\ldots,v_{1}$ every vertex
has at most $k$ preceding neighbors. Now, if $\delta\left(  G\right)  <k,$
then we can add some edge incident to $v_{n},$ and in the ordering
$v_{n},\ldots,v_{1}$ every vertex will still have at most $k$ preceding
neighbors. Hence Lemma B implies that the resulting graph is $k$-degenerate.
This contradicts that $G$ is edge maximal, so we conclude that $\delta\left(
G\right)  =k.$

Next, as $G$ is $k$-degenerate, we have $e\left(  G\right)  \leq kn-\left(
k^{2}+k\right)  /2,$ and Theorem A implies that
\begin{align*}
\mu\left(  G\right)   &  \leq\frac{k-1}{2}+\sqrt{2m-nk+\frac{\left(
k+1\right)  }{4}^{2}}\\
&  \leq\frac{k-1}{2}+\sqrt{2kn-\left(  k^{2}+k\right)  -nk+\frac{\left(
k+1\right)  }{4}^{2}}\\
&  =\mu\left(  S_{n,k}\right)  .
\end{align*}
Note that if $\mu\left(  G\right)  =\mu\left(  S_{n,k}\right)  ,$ then
$e\left(  G\right)  =kn-\left(  k^{2}+k\right)  /2.$

It remains to prove that if $\mu\left(  G\right)  =\mu\left(  S_{n,k}\right)
,$ then $G=S_{n,k}.$ Assume that $G^{\prime}$ is an edge maximal
$k$-degenerate graph of order $n$ containing $G.$ We have
\[
\mu\left(  S_{n,k}\right)  =\mu\left(  G\right)  \leq\mu\left(  G^{\prime
}\right)  \leq\mu\left(  S_{n,k}\right)  ,
\]
and so
\[
\mu\left(  G\right)  =\mu\left(  G^{\prime}\right)  =\mu\left(  S_{n,k}%
\right)  ,
\]
implying in turn that
\[
e\left(  G\right)  =kn-\left(  k^{2}+k\right)  /2=e\left(  G^{\prime}\right)
,
\]
and so $G=G^{\prime}.$

Arguing as above, we see that $G$ is connected and $\delta\left(  G\right)
=k.$ Now, the condition for equality in Theorem A implies that each degree of
$G$ is equal to $k$ or to $n-1.$ Let $x$ be the number of vertices of degree
$n-1$ in $G.$ We have
\[
x\left(  n-1\right)  +\left(  n-x\right)  k=2e\left(  G\right)  =2kn-k^{2}-k,
\]
and so $x=k.$ Clearly, this implies that $G=S_{n,k},$ completing the proof of
Theorem \ref{tA}.
\end{proof}

\bigskip

\begin{proof}
[\textbf{Proof of Theorem \ref{tI}}]Our proof is based on well-known ideas,
used in the same context first in \cite{ElZh00}. Write $r_{k}\left(  B\right)
$ for the $k$'th rowsum of a matrix $B$. Write $d_{1},\ldots,d_{n}$ for the
degrees of $G$, and let $Q$ denote the signless Laplacian of $G.$ Define the
matrix $M$ as%
\[
M=Q^{2}-\left(  \Delta+2\delta-1\right)  Q.
\]
We shall show that the rowsums of $M$ do not exceed $4m-4\left(  n-1\right)
\delta.$ To begin with, recall the well-known fact (see, e.g., \cite{CvSi09})
\[
r_{k}\left(  Q^{2}\right)  =2d_{k}^{2}+2\sum_{\left\{  i,k\right\}  \in
E\left(  G\right)  }d_{i}.
\]
Now, clearly
\[
\sum_{\left\{  i,k\right\}  \in E\left(  G\right)  }d_{i}=\sum_{i\in V\left(
G\right)  }d_{i}-\sum_{\left\{  i,k\right\}  \notin E\left(  G\right)  }%
d_{i}\leq2m-d_{k}-\left(  n-1-d_{k}\right)  \delta,
\]
and so
\begin{align*}
r_{k}\left(  M\right)   &  =2d_{k}^{2}+2\sum_{i\sim k}d_{i}-2\left(
\Delta+2\delta-1\right)  d_{k}\\
&  \leq4m-2d_{k}-2\left(  n-1-d_{k}\right)  \delta+2d_{k}^{2}-2\left(
\Delta+2\delta-1\right)  d_{k}\\
&  =4m-2\left(  n-1\right)  \delta+2d_{k}^{2}-2d_{k}+2\delta d_{k}-2\Delta
d_{k}-4\delta d_{k}+2d_{k}\\
&  =4m-2\left(  n-1\right)  \delta+2d_{k}^{2}-2\left(  \Delta+\delta\right)
d_{k}.
\end{align*}
To eliminate the term $2d_{k}^{2}$ we shall show that for every $k\in V\left(
G\right)  ,$
\[
d_{k}^{2}-\left(  \Delta+\delta\right)  d_{k}\leq-\Delta\delta.
\]
Indeed, the function $f\left(  x\right)  =x^{2}-\left(  \Delta+\delta\right)
x$ is convex and, in view of $\delta\leq d_{k}\leq\Delta,$ we find that
\begin{align*}
d_{k}^{2}-\left(  \Delta+\delta\right)  d_{k}  & =f\left(  d_{k}\right)
\leq\max\left\{  f\left(  \delta\right)  ,f\left(  \Delta\right)  \right\}
\\
& =\max\left\{  \delta^{2}-\left(  \Delta+\delta\right)  \delta,\text{ }%
\Delta^{2}-\left(  \Delta+\delta\right)  \Delta\right\}  \\
& =-\Delta\delta.
\end{align*}
Hence,%
\[
2d_{k}^{2}\leq2\left(  \Delta+\delta\right)  d_{k}-2\Delta\delta,
\]
and, returning back to $r_{k}\left(  M\right)  $, we get%
\begin{align*}
r_{k}\left(  M\right)   &  \leq4m-2\left(  n-1\right)  \delta+2\left(
\Delta+\delta\right)  d_{k}-2\Delta\delta-2\left(  \Delta+\delta\right)
d_{k}\\
&  =4m-2\left(  n-1+\Delta\right)  \delta.
\end{align*}
Now, writing $\mu\left(  M\right)  $ for the largest eigenvalue of $M,$ we see
that%
\[
\mu\left(  M\right)  \leq\max_{k\in V\left(  G\right)  }r_{k}\left(  M\right)
\leq4m-2\left(  n-1+\Delta\right)  \delta.
\]
On the other hand, letting $q=q\left(  G\right)  ,$ it follows that
$q^{2}-\left(  \Delta+2\delta-1\right)  q$ is an eigenvalue of $M$ and so,
\begin{equation}
q^{2}-\left(  \Delta+2\delta-1\right)  q-4m+2\left(  n-1+\Delta\right)
\delta\leq0.\label{qin}%
\end{equation}
Solving this inequality, and using the well-known $q\left(  G\right)
\leq2\Delta,$ we obtain (\ref{bo}).

To prove the condition for equality, first we assume that $G$ is connected. If
equality holds in (\ref{bo}) then Lemma 2.1 of \cite{ElZh00} implies that all
rowsums of $M$ are equal and so
\[
\sum_{\left\{  i,k\right\}  \notin E\left(  G\right)  \text{ and }i\neq
k}d_{i}=\left(  n-1-d_{k}\right)  \delta
\]
for every $k\in V\left(  G\right)  .$ Therefore, if a vertex is not of degree
$n-1,$ then it must be of degree $\delta,$ implying that $G$ is regular or
every vertex is of degree either $\delta$ or $n-1$.

Now suppose that $G$ is not connected and let equality hold in (\ref{qin}).
Assume that $q\left(  G\right)  =q\left(  G_{1}\right)  $ for a component
$G_{1}$ of $G$ and set $n_{1}=v\left(  G_{1}\right)  ,$ $m_{1}=e\left(
G_{1}\right)  ,$ $\Delta_{1}=\Delta\left(  G_{1}\right)  ,$ and $\delta
_{1}=\delta\left(  G_{1}\right)  .$ In view of $\delta_{1}\geq\delta$ and
$\Delta_{1}\leq\Delta,$ inequality (\ref{bo}) and Proposition \ref{pro1} imply
that%
\begin{align*}
q\left(  G\right)   &  \leq\frac{1}{2}\left(  \Delta_{1}+2\delta_{1}%
-1+\sqrt{\left(  \Delta_{1}+2\delta_{1}-1\right)  ^{2}+16m_{1}-8\left(
n_{1}-1+\Delta_{1}\right)  \delta_{1}}\right) \\
&  \leq\frac{1}{2}\left(  \Delta+2\delta-1+\sqrt{\left(  \Delta+2\delta
-1\right)  ^{2}+16m_{1}-8\left(  n_{1}-1+\Delta\right)  \delta}\right) \\
&  \leq\frac{1}{2}\left(  \Delta+2\delta-1+\sqrt{\left(  \Delta+2\delta
-1\right)  ^{2}+16m-8\left(  n-1+\Delta\right)  \delta}\right)  .
\end{align*}
In the above derivation we also use that $16m_{1}-8n_{1}\delta\leq
16m-8n\delta.$ Since the first and the last term of the above chain of
inequalities are equal, it follows that all inequalities are equalities. In
particular, all vertices which do not belong to $G_{1}$ have degree precisely
$\delta.$ Also, either $G_{1}$ is $\delta$-regular or $n_{1}=\Delta+1$ and
every vertex of $G_{1}$ is of degree $\delta$ or $\Delta$. This completes the
proof of the necessity of the condition for equality in (\ref{bo}). The
sufficiency of the condition is simple and we omit the calculations. The proof
of Theorem \ref{tI} is completed.
\end{proof}

\bigskip

\begin{proof}
[\textbf{Proof of Theorem \ref{tQ}}]Our first goal is to prove that $q\left(
G\right)  \leq q\left(  S_{n,k}\right)  .$ As in the proof of Theorem \ref{tA}
we can assume that $\delta=k.$ Using Corollary \ref{cor1} and the fact the
$2m\leq2kn-k^{2}-k,$ we find that
\begin{align*}
q\left(  G\right)   &  \leq\frac{1}{2}\left(  n+2k-2+\sqrt{\left(
n+2k-2\right)  ^{2}+16m-16\left(  n-1\right)  k}\right)  \\
&  \leq\frac{1}{2}\left(  n+2k-2+\sqrt{\left(  n+2k-2\right)  ^{2}%
+16kn-8k^{2}-8k-16\left(  n-1\right)  k}\right)  \\
&  =\frac{1}{2}\left(  n+2k-2+\sqrt{\left(  n+2k-2\right)  ^{2}-8\left(
k^{2}-k\right)  }\right)  \\
&  =q\left(  S_{n,k}\right)  .
\end{align*}
Note that if $q\left(  G\right)  =q\left(  S_{n,k}\right)  ,$ then $e\left(
G\right)  =kn-\left(  k^{2}+k\right)  /2.$

It remains to prove that if $q\left(  G\right)  =q\left(  S_{n,k}\right)  ,$
then $G=S_{n,k}.$ Assume that $G^{\prime}$ is an edge maximal $k$-degenerate
graph of order $n$ containing $G.$ We have
\[
q\left(  S_{n,k}\right)  =q\left(  G\right)  \leq q\left(  G^{\prime}\right)
\leq q\left(  S_{n,k}\right)  ,
\]
and so
\[
q\left(  G\right)  =q\left(  G^{\prime}\right)  =q\left(  S_{n,k}\right)  .
\]
Since
\[
e\left(  G\right)  =kn-\left(  k^{2}+k\right)  /2=e\left(  G^{\prime}\right)
,
\]
it turns out that $G=G^{\prime}.$ Arguing as in the proof of Theorem \ref{tA},
we see that $\delta\left(  G\right)  =k.$ Now, the condition for equality of
Corollary \ref{cor1} implies that each degree of $G$ is equal to $k$ or to
$n-1.$ From here the proof of Theorem \ref{tQ} is completed exactly as the
proof of Theorem \ref{tA}.
\end{proof}

\bigskip

\textbf{Acknowledgement }An idea of Bela Bollob\'{a}s led to the simple proof
of Theorem \ref{tA}.

\end{document}